\documentclass[12pt, draft]{article}
\usepackage{amsmath,amssymb,latexsym,array}
\usepackage[T2A]{fontenc}
\usepackage[cp1251]{inputenc}
\usepackage{longtable}

\renewcommand{\S}{\mathhexbox278}

\DeclareMathOperator{\RRe}{Re}

\newcommand{\Z}{\scriptstyle}
\newcommand{\D}{\displaystyle}

\renewcommand{\le}{\operatorname{\leqslant}}
\renewcommand{\ge}{\operatorname{\geqslant}}

\multlinegap=0em \DeclareFontFamily{T1}{msb}{}
\DeclareFontShape{T1}{msb}{m}{ol}{<5> <6> <7> <8> <9> gen * msbm
<10> <10.95> <12> <14.4> <17.28> <20.74> <24.88> msbm10}{}
\DeclareSymbolFont{AMSb}{T1}{msb}{m}{ol} \multlinegap=0em
\begin{document}
\begin{center}
{\rmfamily\bfseries\normalsize Upper and lower bounds for the
function $\boldsymbol{S(t)}$ on the short intervals\footnote{This
research is supported by the grant of Russian Fund of Fundamental
Researches № 12-01-33080.}}
\end{center}

\begin{center}
\textbf{M.A. Korolev}
\end{center}

\vspace{5mm}

\textbf{Abstract.} We prove under RH the existence of a very large
positive and negative values of the argument of the Riemann zeta
function on a very short intervals.

\vspace{5mm}

In this paper, we study an upper and lower bounds for the function
\[
S(t)\,=\,\pi^{-1}\,\arg{\zeta\bigl(\tfrac{1}{2}+it\bigr)}
\]
on the short intervals. We refer to \cite{Karatsuba_Korolev_2005}
for the definition and for the basic properties of $S(t)$  and to
\cite{Korolev_2005} for the history of the question. Here we mention
only the recent result of R.N.~Boyarinov \cite{Boyarinov_2011}:

\vspace{2mm}

\textsc{Theorem (R.N.~Boyarinov)}. \emph{Let $T>T_{0}>0$ and let
\[
\sqrt{\log\log{T}}\,<\,H\,\le\,(\log{T})(\log\log{T})^{-\,\tfrac{\Z
3}{\Z 2\mathstrut}}.
\]
If the Riemann hypothesis is true then the inequalities
\[
\sup_{T-H\le t\le T+2H}\bigl\{\pm
S(t)\bigr\}\,\ge\,\frac{1}{900}\,\frac{\sqrt{\log{H}\mathstrut}}{\log\log{H}}
\]
hold.}

\vspace{2mm}

In what follows, we prove the similar assertion for the case when
$H$ is essentially smaller than $\sqrt{\log\log{T}\mathstrut}$.
Namely, we prove

\vspace{2mm}

\textsc{Theorem.} \emph{Let $m\ge 2$ be any fixed integer, $T>T_{0}(m)>m$
and let
\[
\frac{2(2m\log\log{T})^{\frac{\Z 1}{\Z
2m\mathstrut}}}{\log\log\log\log{T}}\,\le\,H\,\le\,\sqrt{\log\log{T}}.
\]
If the Riemann hypothesis is true then the inequalities
\[
\sup_{T-H\le t\le T+2H}\bigl\{\pm
S(t)\bigr\}\,\ge\,\frac{1}{50\pi}\,\frac{\sqrt{\log{H}\mathstrut}}{(8m\log\log{H})^{m}}
\]
hold.}

\vspace{2mm}

\textsc{Notations.} We use the following notations:

-- $m\ge 2$ is any fixed integer;

-- $\Phi(u) = \exp{\Bigl(-\,\frac{\D u^{2m}}{\D
2m\mathstrut}\Bigr)}$;

-- $\Lambda(n)$ denotes von Mangoldt's function: $\Lambda(n) =
\log{p}$ if $n = p^{m}$ and $p$ is prime, and $\Lambda(n) = 0$
otherwise;

-- $\widehat{f}$ denotes the Fourier transform of $f$, that is $\D
\widehat{f}(\lambda) =
\int\limits_{-\infty}^{+\infty}f(u)e^{-i\lambda u}du$;

-- $\theta, \theta_{1}, \theta_{2}, \ldots $ denote complex numbers
whose absolute values do not exceed one, different in different
relations.

\vspace{2mm}

We need some auxilliary assertions.

\vspace{2mm}

\textsc{Lemma 1.} \emph{The function $\widehat{\Phi}(\lambda)$
decreases monotonically on the segment $0\le \lambda \le 1$ from the
value $\widehat{\Phi}(0) = \frac{\D 2 \Gamma(1/(2m))}{\D (2m)^{1 -
1/(2m)}}$ to the value $\widehat{\Phi}(1) > \D
\frac{5}{4}\,\Phi\Bigl(\frac{\pi}{4}\Bigr)$. Moreover, there exists
the constant $\lambda_{0} = \lambda_{0}(m)$ such that the inequality
\[
\bigl|\widehat{\Phi}(\lambda)\bigr|\,<\,\frac{5}{\sqrt{m}}\,|\lambda|^{-\beta}\exp{\Bigl(-\,\frac{|\lambda|^{\alpha}}{\alpha}\sin{(\pi\kappa)}\Bigr)}
\]
holds for any real $\lambda$, $|\lambda|>\lambda_{0}$, with}
\[
\alpha\,=\,\frac{2m}{2m-1},\quad \beta = \frac{m-1}{2m-1},\quad
\kappa = \frac{1}{2(2m-1)}.
\]

\vspace{2mm}

\textsc{Proof.} Differentiating the equation
\[
\widehat{\Phi}(\lambda)\,=\,2\int\limits_{-\infty}^{+\infty}\Phi(u)\cos{(\lambda
u)}du
\]
with respect to $\lambda$, we obtain $\widehat{\Phi}'(\lambda) =
-2j(\lambda)$ where
\[
j(\lambda)\,=\,\int\limits_{-\infty}^{+\infty}u\Phi(u)\sin{(\lambda
u)}du.
\]
Suppose that $\lambda > 0$. Using the inequalities $\sin(\lambda
u)\ge \frac{\D 2}{\D \pi}\,\lambda u$ for $0\le u\le \frac{\D
\pi}{\D 2\lambda\mathstrut}$ and $\sin(\lambda u)\ge 0$ for
$\frac{\D \pi}{\D 2\lambda\mathstrut}\le u\le \frac{\D \pi}{\D
\lambda\mathstrut}$ and splitting the integral $j(\lambda)$ into the
sum
\[
\biggl(\;\int\limits_{0}^{\frac{\Z \pi}{\Z
2\lambda\mathstrut}}\;+\;\int\limits_{\frac{\Z \pi}{\Z
2\lambda\mathstrut}}^{\frac{\Z \pi}{\Z
\lambda\mathstrut}}\;+\;\int\limits_{\frac{\Z \pi}{\Z
\lambda\mathstrut}}^{+\infty}\;\;\biggr)\,u\Phi(u)\sin{(\lambda
u)}\,du\,=\,j_{1}+j_{2}+j_{3},
\]
we get
\begin{align*}
& j_{1}\,\ge\,\frac{2\lambda}{\pi}\int\limits_{0}^{\frac{\Z \pi}{\Z
2\lambda\mathstrut}}u^{2}\Phi(u)du\,>\,\frac{2\lambda}{\pi}\,\Phi\Bigl(\frac{\pi}{2\lambda}\Bigr)\int\limits_{0}^{\frac{\Z
\pi}{\Z
2\lambda\mathstrut}}u^{2}du\,=\,\frac{1}{3}\biggl(\frac{\pi}{2\lambda}\biggr)^{2}\Phi\Bigl(\frac{\pi}{2\lambda}\Bigr),\quad
j_{2}\,>\,0;
\end{align*}
\begin{align*}
 & |j_{3}|\,\le\,\int\limits_{\frac{\Z \pi}{\Z
\lambda\mathstrut}}^{+\infty}u\Phi(u)du\,=\,(2m)^{\frac{\Z 1}{\Z
m}-1}\!\!\int\limits_{\frac{\Z 1}{\Z 2m\mathstrut}\left(\frac{\Z
\pi}{\Z \lambda\mathstrut}\right)^{2m}}^{+\infty}w^{\frac{\Z 1}{\Z
m}-1}e^{-w}dw\,<\\
&
<\,\Bigl(\frac{\lambda}{\pi}\Bigr)^{2(m-1)}\!\!\int\limits_{\frac{\Z
1}{\Z 2m\mathstrut}\left(\frac{\Z \pi}{\Z
\lambda\mathstrut}\right)^{2m}}^{+\infty}e^{-w}dw\,=\,\Bigl(\frac{\lambda}{\pi}\Bigr)^{2(m-1)}\Phi\Bigl(\frac{\pi}{\lambda}\Bigr).
\end{align*}
Hence,
\begin{multline*}
j(\lambda)>j_{1}-|j_{3}| >
\frac{1}{3}\Bigl(\frac{\pi}{2\lambda}\Bigr)^{2}\Phi\Bigl(\frac{\pi}{2\lambda}\Bigr)\biggl\{1\,-\,12\biggl(\frac{\lambda}{\pi}\biggr)^{2m}
\exp{\biggl(-\,\frac{1}{2m}\,(1-2^{-2m})\Bigl(\frac{\pi}{\lambda}\Bigr)^{2m}\biggr)}\biggr\}.
\end{multline*}
If $0<\lambda\le 1$ then the expression in the figure brackets is
bounded from below by the value
\[
1\,-\,\frac{12}{\pi^{4}}\,\exp{\biggl(-\,\frac{15}{4}\Bigl(\frac{\pi}{2}\Bigr)^{4}\biggr)}\,>\,1-2\cdot
10^{-10}\,>\,0
\]
uniformly for $m\ge 2$. Therefore, $\Phi'(\lambda)<0$ for
$0<\lambda\le 1$. Thus we prove the first assertion of the lemma.

Next, splitting the expression for $\widehat{\Phi}(1)$ into the sum
\[
2\int\limits_{-\infty}^{+\infty}\Phi(u)\cos{u}\,du\,=\,2\biggl(\;\int\limits_{0}^{\frac{\Z
\pi}{\Z 4\mathstrut}}\;+\;\int\limits_{\frac{\Z \pi}{\Z
4\mathstrut}}^{\frac{\Z \pi}{\Z
2\mathstrut}}\;+\int\limits_{{\frac{\Z \pi}{\Z
2\mathstrut}}}^{+\infty}\;\biggr)\Phi(u)\cos{u}\,du\;=\;2(j_{1}+j_{2}+j_{3})
\]
and using the same arguments as above, we obtain:
\begin{align*}
& j_{1}\,>\,\Phi\Bigl(\frac{\pi}{4}\Bigr)\int\limits_{0}^{\frac{\Z
\pi}{\Z
4\mathstrut}}\cos{u}\,du\,=\,\frac{1}{\sqrt{2}}\Phi\Bigl(\frac{\pi}{4}\Bigr),\quad j_{2}\,>\,0,\\
& |j_{3}|\,\le\,(2m)^{\frac{\Z 1}{\Z
2m\mathstrut}-1}\!\!\int\limits_{\frac{\Z 1}{\Z
2m\mathstrut}\left(\frac{\Z \pi}{\Z
2\mathstrut}\right)^{2m}}^{+\infty}w^{\frac{\Z 1}{\Z
2m\mathstrut}-1}e^{-w}dw\,\le\,\biggl(\frac{2}{\pi}\biggr)^{2m}\Phi\Bigl(\frac{\pi}{2}\Bigr),
\end{align*}
and hence
\begin{multline*}
\widehat{\Phi}(1)\,>\,j_{1}-|j_{3}|\,>\,2\biggl\{\frac{1}{\sqrt{2}}\,\Phi\Bigl(\frac{\pi}{4}\Bigr)\,-\,\biggl(\frac{2}{\pi}\biggr)^{4}\Phi\Bigl(\frac{\pi}{2}\Bigr)\biggr\}\,=\\
=\,\sqrt{2}\,\Phi\Bigl(\frac{\pi}{4}\Bigr)\Biggl\{1\;-\;\sqrt{2}\biggl(\frac{2}{\pi}\biggr)^{4}\,\frac{\Phi\Bigl(\frac{\D
\pi}{\D 2\mathstrut}\Bigr)}{\Phi\Bigl(\frac{\D \pi}{\D
4\mathstrut}\Bigr)}\Biggr\}.
\end{multline*}
One can note that
\[
\frac{\Phi\Bigl(\frac{\D \pi}{\D
2\mathstrut}\Bigr)}{\Phi\Bigl(\frac{\D \pi}{\D
4\mathstrut}\Bigr)}\,=\,\exp{\biggl\{-\,\frac{1}{2m}\,(1-2^{-2m})\biggl(\frac{\pi}{2}\biggr)^{2m}\biggr\}}\,\le\,\exp{\biggl(-\,\frac{15}{64}\biggl(\frac{\pi}{2}\biggr)^{4}\biggr)}
\]
for any $m\ge 2$. Thus we arrive at the desired bound for
$\widehat{\Phi}(1)$.

Finally, the last assertion of the lemma follows from the formula
\begin{multline*}
\int\limits_{-\infty}^{+\infty}\exp{\biggl(-\,\frac{u^{2m}}{2m}+i\lambda
u\biggr)}du\,=\\
=\,\frac{4\sqrt{\pi\kappa}}{\lambda^{\beta\mathstrut}}\,\exp{\biggl(-\,\frac{\lambda^{\alpha}}{\alpha}\,\sin{(\pi\kappa)}\biggr)}\,
\biggl\{\cos{\biggl(\frac{\lambda^{\alpha}}{\alpha}\,\cos{(\pi\kappa)}\biggr)}\,+\,O\bigl(\lambda^{-\,\alpha}\bigr)\biggr\},
\end{multline*}
where $\lambda\to +\infty$ and $\alpha, \beta, \kappa$ are defined
as above (see \cite[\S 7.1]{Fedoryuk_1977}).

\vspace{2mm}

\textsc{Lemma 2.} \emph{The following inequalities hold true:}
\begin{equation*}
|S(t)|\,\le\,
\begin{cases}
1, & \emph{if}\;\;|t|\le 280, \\
1.05\log{|t|}, & \emph{if}\;\; |t|>280.
\end{cases}
\end{equation*}

\vspace{2mm}

\textsc{Proof.} The first inequality follows from the data of Table
1 from \cite{Lehman_1970} and the second one follows from the
classical estimate of R.J.~Backlund \cite{Backlund_1916}:
\begin{multline*}
|S(t)|\,<\,0.1361\log{|t|}\,+\,0.4422\log\log{|t|}\,+\,4.3451\,\le\\
\le\,(\log{|t|})\biggl(0.1361\,+\,0.4422\,\frac{\log\log{280}}{\log{280}}\,+\,\frac{4.3451}{\log{280}}\biggr)\,<\,1.05\log{|t|}.
\end{multline*}
The lemma is proved.

\vspace{2mm}

Let $\tau > 1$ and $f(u) = \Phi(\tau u)$.

\vspace{2mm}

\textsc{Lemma 3.} \emph{If the Riemann hypothesis is true then the formula
\[
\int\limits_{-\infty}^{+\infty}f(u)S(t+u)du\,=\,-\,\frac{1}{\pi}\sum\limits_{n
=
2}^{+\infty}\frac{\Lambda(n)}{\sqrt{n}}\,\frac{\sin{(t\log{n})}}{\log{n}}\,\widehat{f}(\log{n})\,-\,2\int\limits_{0}^{1/2}f(-t-iu)du
\]
holds.}

\vspace{2mm}

\textsc{Proof.} The proof of this assertion repeats word\,-by\,-word
(with the minor changes) the proof of Theorem 3 from
\cite{Karatsuba_Korolev_2005}. The difference is that we need to use
the explicit formula for $f(z)$ instead of the inequality $|f(z)|\le
c(|z|+1)^{-(1+\alpha)}$.

\vspace{2mm}

\textsc{Lemma 4.} \emph{Let $y>y_{0}>0$, $\mu, \nu\ge 0, k\ge 1, k =
\mu + \nu$, and let $p_{1},\ldots, p_{\mu}$, $q_{1},\ldots, q_{\nu}$
range over the primes in the interval $(1,y]$ that satisfy the
condition $p_{1}\ldots p_{\mu}\ne q_{1}\ldots q_{\nu}$. If $a(p)$ is
a sequence of complex numbers in which $|a(p)|\le \delta$ for any
prime $p\le y$, then the integral
\[
I\,=\,\int\limits_{T}^{T+H}\sum\limits_{\substack{p_{1},\ldots,
p_{\mu} \\ q_{1},\ldots, q_{\nu}}}\frac{a(p_{1})\ldots
a(p_{\mu})\,\overline{a}(q_{1})\ldots
\overline{a}(q_{\nu})}{\sqrt{p_{1}\ldots
q_{\mu}}}\,\biggl(\frac{p_{1}\ldots p_{\mu}}{q_{1}\ldots
q_{\nu}}\biggr)^{\! it}dt
\]
satisfies the estimate} $|I|\le (\delta^{2}y^{3})^{k}$.

\vspace{2mm}

\textsc{Proof.} It is lemma 2 from \cite[\S
2.1]{Karatsuba_Korolev_2006}.

\vspace{2mm}

\textsc{Lemma 5.} \emph{Let $k\ge 1$ be an integer, let $M>0$ and
let a real function $W(t)$ satisfies the inequalities
\[
\int\limits_{T}^{T+H}W^{2k}(t)\,dt\,\ge\,HM^{2k},\quad
\biggl|\int\limits_{T}^{T+H}W^{2k+1}(t)\,dt\biggr|\,\le\,\frac{1}{2}\,HM^{2k+1}.
\]
Then}
\[
\sup_{T\le t\le T+H}\bigl\{\pm\,W(t)\bigr\}\,\ge\,\frac{M}{2}.
\]

\vspace{2mm}

\textsc{Proof.} It is a slight modification of lemma 4 from the
paper of K.\,-M.~Tsang \cite{Tsang_1986}.

\vspace{2mm}

\textsc{Proof of the theorem.} Let $\tau = 2\log\log{H}$ and suppose
that $T\le t\le T+H$. Applying lemma 3 to the function $f(u) =
\Phi(\tau u)$ we obtain
\[
\tau\int\limits_{-\infty}^{+\infty}\Phi(\tau
u)S(t+u)du\,=\,-\,\frac{1}{\pi}\,V(t)\,+\,R_{1}(t),
\]
where
\[
V(t)=\sum\limits_{n =
2}^{+\infty}\frac{\Lambda(n)}{\sqrt{n}}\,\frac{\sin{(t\log{n})}}{\log{n}}\,\widehat{\Phi}\Bigl(\frac{\log{n}}{\tau}\Bigr),\quad
R_{1}(t) = 2\int\limits_{0}^{1/2}\Phi\bigl(\tau(t+iu)\bigr)\,du.
\]
Since
\[
\RRe\,(t+iu)^{2m}\,=\,t^{2m}\biggl(1\,+\,\sum\limits_{\nu =
1}^{m}(-1)^{\nu}\binom{2m}{2\nu}\biggl(\frac{u}{t}\biggr)^{2m}\biggr)
\]
and since the absolute value of the last sum is bounded from above
by the values
\[
\sum\limits_{\nu =
1}^{2m}\binom{2m}{\nu}(2t)^{-\,2\nu}\,=\,\biggl(1\,+\,\frac{1}{2t}\biggr)^{2m}-1\,\le\,\frac{1}{2t}\biggl(1+\frac{1}{2t}\biggr)^{2m-1}<\frac{e}{2t},
\]
we find that
\begin{align*}
&
\bigl|\Phi\bigr(\tau(t+iu)\bigl)\bigr|\,=\,\exp{\biggl\{-\,\frac{\tau^{2m}}{2m}\,\RRe(t+iu)^{2m}\biggr\}}\,<\,\exp{\biggl\{-\,\frac{(\tau
t)^{2m}}{8m}\biggr\}},\\
& |R_{1}(t)|\,\le\,2\cdot\frac{1}{2}\,\exp{\biggl\{-\,\frac{(\tau
t)^{2m}}{8m}\biggr\}}.
\end{align*}
Further, let $X = \exp{\bigl((4m\tau)^{2m}\bigr)}$ and let
$R_{2}(t)$ be the contribution to the sum $V(t)$ from the terms with
$n>X$. Since $\frac{\D \log{n}}{\D \tau} > \lambda_{0}$ for any
$n>X$ ($\lambda_{0}$ is defined in lemma 1) then lemma 1 implies
that
\begin{equation}\label{Lab_01}
\biggl|\widehat{\Phi}\Bigl(\frac{\log{n}}{\tau}\Bigr)\biggr|\,\le\,\frac{5}{\sqrt{m}}\biggl(\frac{\tau}{\log{n}}\biggr)^{\beta}\,
\exp{\biggl\{-\,\frac{\sin{(\pi\kappa)}}{\alpha}\biggl(\frac{\log{n}}{\tau}\biggr)^{\alpha}\biggr\}},
\end{equation}
where the values $\alpha, \beta$ and $\kappa$ are defined above. By
the inequality $\frac{\D \sin{(\pi\kappa)}}{\D \alpha} \ge \frac{\D
2}{\D \pi}\,\frac{\D\pi\kappa}{\D\alpha} = \frac{\D 1}{\D
2m\mathstrut}$ we conclude from (\ref{Lab_01}) that
\[
|R_{2}(t)|\,<\,\frac{5}{\sqrt{m}}\biggl(\frac{\tau}{\log{X}}\biggr)^{\beta}\sum\limits_{n>X}\frac{1}{\sqrt{n}}\,
\exp{\biggl\{-\,\frac{1}{2m}\biggl(\frac{\log{n}}{\tau}\biggr)^{\alpha}\biggr\}}.
\]
In order to estimate the last sum,we take $r_{0} = \Bigl[\frac{\D
\log{X}}{\D \tau}\Bigr] = \bigl[(4m)^{2m}\tau^{2m-1}\bigr]$ and
split the domain of summation into the segments of the type
$e^{r\tau}<n\le e^{(r+1)\tau}$, $r = r_{0}, r_{0}+1,\ldots$. The sum
over such segment is bounded by the values
\begin{multline*}
\exp{\biggl\{-\,\frac{r^{\alpha}}{2m}\biggr\}}\sum\limits_{e^{r\tau}<n\le
e^{(r+1)\tau}}\frac{1}{\sqrt{n}} <
2\exp{\biggl\{\frac{(r+1)\tau}{2}\,-\,\frac{r^{\alpha}}{2m}\biggr\}}\,<\\
<\,\exp{\biggl\{r\tau\,-\,\frac{r^{\alpha}}{2m}\biggr\}}\,<\,\exp{\biggl\{-\,\frac{r^{\alpha}}{4m}\biggr\}}.
\end{multline*}
Hence,
\begin{multline*}
|R_{2}(t)|\,<\,\frac{5}{\sqrt{m}}\biggl(\frac{\tau}{\log{X}}\biggr)^{\beta}\int\limits_{(4m\tau)^{2m-1}}^{+\infty}\exp{\biggl\{-\,\frac{u^{\alpha}}{2m}\biggr\}}\,du\,=\\
=\,\frac{5}{\sqrt{m}}\biggl(\frac{\tau}{(4m\tau)^{2m}}\biggr)^{\beta}\,\frac{(2m)^{\frac{\Z
1}{\Z \alpha}}}{\alpha}\int\limits_{\frac{\Z 1}{\Z
2m\mathstrut}(4m\tau)^{2m}}^{+\infty}w^{\frac{\Z 1}{\Z
\alpha}-1}e^{-\,w}dw\,<\\
<\,\frac{5}{\sqrt{m}}\biggl(\frac{\tau}{(4m\tau)^{2m}}\biggr)^{\beta}\,\frac{1}{2\tau}\,\exp{\biggl\{-\,\frac{(4m\tau)^{2m}}{2m}\biggr\}}\,=\\
=\,\frac{5(4m)^{\frac{\Z 1}{\Z
2(2m-1)\mathstrut}}}{(4m\tau)^{m}}\,\exp{\biggl\{-\,\frac{(4m\tau)^{2m}}{2m}\biggr\}}\,
\le\,\frac{5\sqrt{2}}{(4m\tau)^{m}}\,\exp{\biggl\{-\,\frac{(4m\tau)^{2m}}{2m}\biggr\}}.
\end{multline*}
Finally, let $R_{3}(t)$ be the sum over $n = p^{\nu}\le X$ with
$\nu\ge 2$. Then the obvious estimate
\[
\biggl|\widehat{\Phi}\Bigl(\frac{\log{n}}{\tau}\Bigr)\biggr|\,=\,\biggl|\int\limits_{-\infty}^{+\infty}\Phi(u)n^{-\,\frac{\Z
iu}{\Z \tau}}\,du\biggr|\,\le\,\widehat{\Phi}(0)
\]
yields
\begin{multline*}
|R_{3}(t)|\,\le\,\widehat{\Phi}(0)\sum\limits_{\nu =
2}^{+\infty}\sum\limits_{p\le
X^{1/\nu}}\frac{p^{-\,\nu/2}}{\nu}\,\le
\,\widehat{\Phi}(0)\biggl(\;\frac{1}{2}\sum\limits_{p\le
\sqrt{X}}\frac{1}{p}\,+\,\frac{1}{3}\sum\limits_{\nu =
3}^{+\infty}\frac{p^{-\,\nu/2}}{\nu}\;\biggr)\,<\\
<\widehat{\Phi}(0)\biggl(\;\frac{1}{2}\log\log{X}\,+\,c\biggr)\,=\,\widehat{\Phi}(0)\bigl(m\log{(4m\tau)}+c\bigr),
\end{multline*}
where $c>0$ is a sufficiently large absolute constant.

Thus we get
\begin{equation*}
\tau\int\limits_{-\infty}^{+\infty}\Phi(\tau
u)S(t+u)\,du\,=\,-\,\frac{1}{\pi}\,W(t)\,+\,\theta_{1}Q_{1},
\end{equation*}
where
\begin{align*}
& W(t)\,=\,\sum\limits_{p\le
x}\frac{a(p)}{\sqrt{p}}\,\sin{(t\log{p})},\quad
a(p)\,=\,\Phi\Bigl(\frac{\log{p}}{\tau}\Bigr),\\
& Q_{1}\,=\,\exp{\biggl\{-\,\frac{(\tau
t)^{2m}}{8m}\biggr\}}\,+\,\frac{5\sqrt{2}}{(4m\tau)^{m}}\,\exp{\biggl\{-\,\frac{(4m\tau)^{2m}}{2m}\biggr\}}+\widehat{\Phi}(0)\bigl(m\log{(4m\tau)}+c\bigr)\,<\\
& <\,\widehat{\Phi}(0)\bigl(m\log{(4m\tau)}+2c\bigr).
\end{align*}

Now let us consider the integrals
\[
I_{1}\,=\,\int\limits_{H}^{+\infty}\Phi(\tau u)S(t+u)\,du, \quad
I_{2}\,=\,\int\limits_{-\infty}^{-H}\Phi(\tau u)S(t+u)\,du.
\]
Splitting $I_{1}$ into the sum
\[
I_{1}\,=\,\biggl(\;\int\limits_{H}^{T}+\;\int\limits_{T}^{+\infty}\biggr)\Phi(\tau
u)S(t+u)\,du\,=\,I_{1}^{(1)}\,+\,I_{1}^{(2)}
\]
and applying lemma 2 we obtain:
\begin{align*}
& |I_{1}^{(1)}|\,\le\,1.05\int\limits_{H}^{T}\Phi(\tau
u)\log{(t+u)}\,du\,<\,\frac{1.1\log{T}}{\tau}\int\limits_{H\tau}^{+\infty}\Phi(v)\,dv\,\le\,\frac{1.1}{\tau}\,\frac{\Phi(H\tau
)\log{T}}{(H\tau)^{2m-1\mathstrut}},\\
& |I_{1}^{(2)}|\,\le\,1.05\int\limits_{T}^{+\infty}\Phi(\tau
u)\log(t+u)\,du\,<\,\frac{1.1}{\tau}\int\limits_{T\tau}^{+\infty}\Phi(v)\log{v}\,dv\,\le\,\frac{1.1}{\tau}\,\frac{2m\Phi(T\tau
)\log{(T\tau)}}{(T\tau)^{2m-1\mathstrut}},
\end{align*}
and hence
\[
|I_{1}|\,<\,\frac{1.2}{\tau}\,\frac{\Phi(H\tau)\log{T}}{(H\tau)^{2m-1}}.
\]

Next, we split the integral $I_{2}$ into the sum
\begin{multline*}
\int\limits_{H}^{+\infty}\Phi(\tau
u)S(t-u)\,du\,=\,\biggl(\;\int\limits_{H}^{t-10^{2}}\,+\,\int\limits_{t-10^{2}}^{t+10^{2}}\,+\int\limits_{t+10^{2}}^{+\infty}\;\biggr)\Phi(\tau
u)S(t-u)\,du\,=\\
=\, I_{2}^{(1)}\,+\,I_{2}^{(2)}\,+\,I_{2}^{(3)}.
\end{multline*}
Applying the first inequality of lemma 2 to $I_{2}^{(2)}$ and the
second one to the estimation of $I_{2}^{(1)}$ and $I_{2}^{(3)}$ we
obtain
\[
|I_{2}|\,<\,\frac{1.2}{\tau}\,\frac{\Phi(H\tau)\log{T}}{(H\tau)^{2m-1}},\quad
|I_{1}|+|I_{2}|\,<\,\frac{2.4}{\tau}\,\frac{\Phi(H\tau)\log{T}}{(H\tau)^{2m-1}}.
\]
Finally, we have
\[
j(t)\,=\,\tau\int\limits_{-H}^{H}\Phi(\tau
u)S(t+u)\,du\,=\,-\,\frac{1}{\pi}\,W(t)\,+\,\theta\,Q_{2},
\]
where
\[
Q_{2}\,=\,\widehat{\Phi}(0)\bigl(m\log{(4m\tau)}\,+\,2c\bigr)\,+\,\frac{2.4}{\tau}\,\frac{\Phi(H\tau)\log{T}}{(H\tau)^{2m-1}}.
\]
Since $H\tau > (2m\log\log{T})^{\frac{\Z 1}{\Z 2m\mathstrut}}$ for
$H$ and $m$ under considering, we find that
\[
\Phi(H\tau)\,<\,(\log{T})^{-1},\quad
Q_{2}\,<\,2\widehat{\Phi}(0)m\log{(4m\tau)}.
\]
Now let us take $k = \Bigl[\frac{\D \log{H}}{\D
5\log{X}\mathstrut}\Bigr]$ and define the integrals $I(k)$ and
$J(k)$ by the following relations:
\[
I(k)\,=\,\int\limits_{T}^{T+H}W^{2k}(t)\,dt,\quad
J(k)\,=\,\int\limits_{T}^{T+H}W^{2k+1}(t)\,dt.
\]
Writing $W(t)$ as
\[
\frac{1}{2i}\,\bigl(U(t)\,-\,\overline{U}(t)\bigr),\quad
U(t)\,=\,\sum\limits_{p\le X}\frac{a(p)}{\sqrt{p}}\,p^{it},
\]
we find that
\[
I(k)\,=\,(2i)^{-2k}\sum\limits_{\nu =
0}^{2k}(-1)^{\nu}\binom{2k}{\nu}j_{\nu},\quad
j_{\nu}\,=\,\int\limits_{T}^{T+H}U^{\,\nu}(t)\overline{U}^{\,2k-\nu}(t)dt.
\]
The application of lemma 4 with $\delta = \widehat{\Phi}(0)$ to the
case $\nu\ne k$ yields:
\[
|j_{\nu}|\,<\,\Bigl(\widehat{\Phi}(0)X^{\frac{\Z 3}{\Z
2\mathstrut}}\Bigr)^{2k}<X^{4k}\le H^{\frac{\Z 4}{\Z 5\mathstrut}}.
\]
The same estimate is valid for the contribution to $j_{k}$ of the
terms under the condition $p_{1}\ldots p_{k}\ne q_{1}\ldots q_{k}$.
Hence,
\[
I(k)\,=\,2^{-2k}\binom{2k}{k}H\mathfrak{S}_{k}\,+\,\theta
H^{\frac{\Z 4}{\Z 5\mathstrut}},
\]
where
\[
\mathfrak{S}_{k}\,=\,\sum\limits_{\substack{p_{1}\ldots p_{k}=
q_{1}\ldots q_{k} \\ p_{1},\ldots, q_{k}\le
X}}\frac{a^{2}(p_{1})\ldots a^{2}(p_{k})}{p_{1}\ldots p_{k}}.
\]
In order to estimate the sum $\mathfrak{S}_{k}$ from the below, we
truncate the sum by replacing the upper bound for $p_{1},\ldots,
q_{k}$ by the value $Y = e^{\tau} = (\log{H})^{2}<X$. Then it
follows from lemma 1 that the inequalities
\[
a(p)\,=\,\widehat{\Phi}\Bigl(\frac{\log{p}}{\tau}\Bigr)\,\ge\,\widehat{\Phi}(1)
\]
hold for any $p\le Y$. Further, if we retain the terms in the
truncated sum that correspond to the tuples $(p_{1},\ldots, p_{k})$
involving no repetitions and using the fact that the number of
solutions $(q_{1},\ldots, q_{k})$ of the equation $p_{1}\ldots
p_{k}= q_{1}\ldots q_{k}$ is equal to $k!$, we obtain that
\[
\mathfrak{S}_{k}\,=\,\bigl(\widehat{\Phi}(1)\bigr)^{2k}\sum\limits_{\substack{p_{1}\ldots
p_{k}= q_{1}\ldots q_{k} \\ p_{1},\ldots, q_{k}\le Y}}(p_{1}\ldots
p_{k})^{-\,1}\,\ge\,k!\,\bigl(\widehat{\Phi}(1)\bigr)^{2k}\sum\limits_{\substack{p_{1},\ldots,\,
p_{k}\le Y \\ p_{1},\ldots, p_{k}\;\text{are distinct}}}(p_{1}\ldots
p_{k})^{-\,1}.
\]
Applying the arguments from \cite{Korolev_2005} (the estimate of the
sum $\Sigma$) to the estimation of the last sum, we find that
\begin{align*}
&
\mathfrak{S}_{k}\,\ge\,k!\,\bigl(\widehat{\Phi}(1)\bigr)^{2k}\biggl(\,\sum\limits_{2k\log{k}<p\le
Y}\frac{1}{p}\biggr)^{k}\,\ge\,k!\,\bigl(\widehat{\Phi}(1)\bigr)^{2k}\biggl(\,\sum\limits_{\sqrt{Y}<p\le
Y}\frac{1}{p}\biggr)^{k}\,>\,k!\,\biggl(\frac{4}{5}\,\widehat{\Phi}(1)\biggr)^{2k},\\
&
I(k)\,>\,\frac{(2k)!}{k!}\,\frac{H}{2^{2k\mathstrut}}\,\biggl(\frac{4}{5}\,\widehat{\Phi}(1)\biggr)^{2k}\,-\,H^{\frac{\Z
4}{\Z 5\mathstrut}}\,
>\,\frac{e}{2}\biggl(\frac{4k}{e}\biggr)^{k}\biggl(\frac{\widehat{\Phi}(1)}{5}\biggr)^{2k}\,-\,H^{\frac{\Z
4}{\Z 5\mathstrut}}\,>\, HM^{2k},
\end{align*}
where
\[
M\,=\,\frac{2}{5}\,\widehat{\Phi}(1)\,\sqrt{\frac{k}{e}}\,>\,2.
\]
Finally, lemma 4 yields:
\[
|J(k)|\,<\,\Bigl(\widehat{\Phi}(0)X^{\frac{\Z 3}{\Z
2\mathstrut}}\Bigr)^{2k+1}\,<\,X^{4k}\le H^{\frac{\Z 4}{\Z
5\mathstrut}}\,<\,\frac{1}{2}\,HM^{2k+1}.
\]
Now it follows from lemma 5 that there exist the values $t_{0}$ and
$t_{1}$ such that $T\le t_{0}, t_{1}\le T+H$ such that $W(t_{0})<
-0.5M$ and $W(t_{1})>0.5M$. Thus we have
\begin{align*}
&
j(t_{0})\,>\,-\,\frac{1}{\pi}\,W(t_{0})\,-\,Q_{2}\,>\,\frac{M}{2\pi}\,-\,2\widehat{\Phi}(0)m\log{(4m\tau)},\\
&
j(t_{1})\,<\,-\,\frac{1}{\pi}\,W(t_{1})\,+\,Q_{2}\,<\,-\,\frac{M}{2\pi}\,+\,2\widehat{\Phi}(0)m\log{(4m\tau)}.
\end{align*}
Setting $M_{j} = \D \sup_{|u|\le H}(-1)^{j}S(t_{j}+u)$ for $j =
0,1$, we obviously have
\[
j(t_{0})\,<\,M_{0}\tau\int\limits_{-H}^{H}\Phi(\tau
u)\,du\,<\,M_{0}\widehat{\Phi}(0),\quad
j(t_{1})\,>\,M_{1}\widehat{\Phi}(0)
\]
and therefore
\[
(-1)^{j}M_{j}\,>\,\mu,\quad
\mu\,=\,\widehat{\Phi}(0)\,\frac{M}{2\pi}\,-\,2m\log{(4m\tau)}.
\]
Finally, applying lemma 1 together with the inequality
\[
\Gamma\biggl(\frac{1}{2m}\biggr)\,=\,2m\Gamma\biggl(1+\frac{1}{2m}\biggr)\,\le\,2m,
\]
we obtain
\begin{multline*}
\mu\,\ge\,\frac{1}{5\pi}\,\frac{\widehat{\Phi}(1)}{\widehat{\Phi}(0)}\,\sqrt{\frac{k}{e}}\,-\,2m\log{(4m\tau)}\,>\\
>\,\frac{1}{8\pi}\,\frac{\Phi\Bigl(\frac{\D
\pi}{\D 4\mathstrut}\Bigr)}{\Gamma\Bigl(\frac{\D 1}{\D
2m\mathstrut}\Bigr)}\,(2m)^{1\,-\,\frac{\Z 1}{\Z
2m\mathstrut}}\sqrt{\frac{1}{e}\biggl(\frac{\log{H}}{5(4m\tau)^{2m}}\,-\,1\biggr)}\,-\,2m\log{(4m\tau)}\,>\\
>\,\frac{1}{32\pi}\,\frac{\Phi\Bigl(\frac{\D
\pi}{\D 4\mathstrut}\Bigr)}{2m}\,(2m)^{1\,-\,\frac{\Z 1}{\Z
2m\mathstrut}}\,\frac{\sqrt{\log{H}\mathstrut}}{(4m\tau)^{m}}\,\ge\\
\ge\,\frac{1}{32\pi\sqrt{2}}\,\exp{\biggl\{-\,\frac{1}{4}\Bigl(\frac{\pi}{4}\Bigr)^{4}\biggr\}}\,\frac{\sqrt{\log{H}\mathstrut}}{(8m\log\log{H})^{m}}\,>
\,\frac{1}{50\pi}\,\frac{\sqrt{\log{H}\mathstrut}}{(8m\log\log{H})^{m}}.
\end{multline*}
The theorem is proved.

\vspace{2mm}

\textsc{Remark.} The assertion of the theorem can be generalized to
the case when $m$ grows with $T$.

\renewcommand{\refname}{

\end{document}